\documentclass[reqno]{amsart}

\usepackage{microtype, hyperref, cleveref}
\usepackage{amsmath, amssymb}
\usepackage{graphicx}

\newtheorem{con}{Conjecture}
\newtheorem{crl}[con]{Corollary}

\newtheorem{lem}[con]{Lemma}
\newtheorem{thrm}[con]{Theorem}
\newtheorem{rem}[con]{Remark}

\makeatletter
\newcommand\numbereq{%
  \ifmeasuring@\else
    \refstepcounter{equation}%
  \fi
  \tag{\theequation}%
}
\makeatother

\newcommand{\ind}{\operatorname{ind}}
\newcommand{\Spinc}{\operatorname{Spin}^\mathrm{c}}
\newcommand{\Jcal}{\mathcal{J}}
\newcommand{\sfr}{\mathfrak{s}}
\newcommand{\osfr}{\overline{\mathfrak{s}}}
\newcommand{\Nbb}{\mathbb{N}}
\newcommand{\Zbb}{\mathbb{Z}}
\newcommand{\CPtwo}{\mathbb{C}P^2}
\newcommand{\oCPtwo}{\overline{\mathbb{C}P^2}}
\newcommand{\Hbb}{\mathbb{H}}
\renewcommand{\c}{\textrm{c}}

\usepackage[backend=biber,style=alphabetic]{biblatex}
\title{Almost complex structures on triple blow-ups}
\author{Samuel Adrian Antz}
\date{\today}
\address{School of Mathematical Sciences, Dalian University of Technology, P.R. China}
\email{samueladrianantz@mail.dlut.edu.cn}

\begin{document}

\begin{abstract}
We show that for certain smooth 4-manifolds, almost complex structures on its triple blow-up induce all its spin$^\c$ structures with a non-negative-dimensional Seiberg--Witten moduli space. Since this includes all basic classes, our result provides an approximation to the simple type conjecture. Moreover, we expand a result by Dieter Kotschick showing that every finitely presented group is the fundamental group of a closed almost complex 4-manifold. We show that it can furthermore be chosen to have every even integer as the virtual dimension of the Seiberg--Witten moduli space of a spin$^\c$ structure.
\end{abstract}

\maketitle

\section*{Introduction}

Seiberg--Witten theory is the application of the eponymous gauge theory to smooth $4$-manifolds $M$, as described by Edward Witten in \cite{witten94}. It first requires a spin$^\c$ structure $\sfr$ to formulate the Seiberg--Witten equations. It then considers the Seiberg--Witten moduli space of all solutions up to gauge transformation. Such a moduli space is associated a virtual dimension $d(\sfr)\in\Zbb$, dictated purely by the topology of $M$ and the spin$^\c$ structure $\sfr$. Almost complex structures are naturally related to this construction: every almost complex structure $J$ induces a spin$^\c$ structure $\sfr_J$ and furthermore $d(\sfr_J)=0$ (see \parencites[p. 459]{nicolaescu00}[p. 390]{scorpan05}). In reverse, every such spin$^\c$ structure is induced by an almost complex structure, giving a surjection.

Our main result provides generalized surjections between spin$^\c$ structures with inequality conditions for the virtual dimension. It is obtained by combining connected sums with $\oCPtwo$ and $\CPtwo$ with Gauss' Eureka theorem (see \cite[Thrm. 6.5.5]{anglin95}) from number theory (assuming $0\in\Nbb$). In particular, it extends the previous statement by showing that every spin$^\c$ structure $\sfr$ with $d(\sfr)\geq 0$ is induced by an almost complex structure on the triple blow-up, meaning a connected sum of the $4$-manifold with $3\oCPtwo$. Its main significance comes from its uniformity. Namely, no matter the manifold or the virtual dimension, three blow-ups are always sufficient to guarantee this and are optimal in general.

A direct consequence of this result concerns the (Seiberg--Witten) simple type conjecture. It claims that for every oriented closed smooth $4$-manifold with $-b_1+b_2^+$ odd and $b_2^+\geq 2$, every (Seiberg--Witten) basic class $\sfr$ (a spin$^\c$ structure with a non-trivial Seiberg--Witten invariant) has $d(\sfr)=0$ or equivalently is induced by an almost complex structure \cite[p. 27]{witten94} (see also \parencites[p. 150]{nicolaescu00}[p. 402]{scorpan05}). Since a basic class $\sfr$ requires $d(\sfr)\geq 0$ by definition, our result can be applied.

Moreover, we also consider the case when the orientation of two blow-ups is changed, meaning a connected sum of the $4$-manifold with $2\CPtwo\#\oCPtwo$. We show that this removes all inequality conditions for the virtual dimension: every spin$^\c$ structure and every even number can be combined to a spin$^\c$ structure on this connected sum. In partiular, every spin$^\c$ structure comes from an almost complex structure on this connected sum. We use this to extend a result by Dieter Kotschick in \cite{kotschick92}, stating that every finitely presented group is the fundamental group of a closed almost complex $4$-manifold. His result was obtained analogous by combining connected sums with $\oCPtwo$ and $\CPtwo$ with Lagrange's four-squares theorem (see \cite[Thrm. 1.8.4]{anglin95}) from number theory. Both Gauss' and Lagrange's theorem are special cases of Fermat's polygonal number theorem (see \cite[Thrm. 6.6.3]{anglin95}), giving further common ground.

The paper is organized as follows: Section \ref{section:results} will state the main result concerning the generalized surjections and their implications. Section \ref{section:setting} will review their setting and Section \ref{section:proofs} will provide their proofs. Section \ref{section:discussion} will discuss aspects like potential improvement, relation to other works and relevant details concerning the simple type conjecture.

\section{Results}
\label{section:results}

Throughout this paper, all manifolds are assumed to be connected. Let $M\# N$ be the connected sum of smooth $4$-manifolds $M$ and $N$. It is obtained by removing an embedded $4$-disc from both, glueing the boundary spheres together and uniquely smoothing it (see \cite[p. 12]{scorpan05}). Collapsing one of the manifolds to a point and smoothing the map by homotopy yields canonical projections $M\# N\twoheadrightarrow M,N$. 

Let $M$ be an oriented closed smooth $4$-manifold with $(-b_1+b_2^+)(M)$ odd. Let $\Spinc(M)$ denote the set of spin$^\c$ structures of $M$. The main result of the paper is the following:
\begin{thrm}
For every $n\in\Nbb$, the canonical projections $p\colon M\#(n+3)\oCPtwo\twoheadrightarrow M$ and $q\colon M\#(2n+4)\CPtwo\twoheadrightarrow M$ induce surjections:
\begin{equation}
(p^*,d)\colon
\Spinc\left(M\#(n+3)\oCPtwo\right)\twoheadrightarrow\left\{(\sfr,2m)\in\Spinc(M)\times 2\Zbb|d(\sfr)\geq 2m\right\};
\label{eq:main_theorem_1}
\end{equation}
\begin{equation}
\begin{aligned}
(q^*,d)\colon
&\Spinc\left(M\#(2n+4)\CPtwo\right) \\
\twoheadrightarrow &\left\{(\sfr,2m)\in\Spinc(M)\times 2\Zbb|d(\sfr)\leq 2(m+n+2)\right\}.
\end{aligned}
\label{eq:main_theorem_2}
\end{equation}
\label{thrm:main_theorem}
\end{thrm}
Induced means that the spin$^\c$ structure on the connected sum is restricted onto one of the manifolds with an embedded $4$-disc removed and then uniquely extended, using that the $3$-sphere has a unique spin$^\c$ structure.

Let $\Jcal(M)$ denote the set of almost complex structures of $M$. For every $m\in\Zbb$, let $\Spinc_m(M)$, $\Spinc_{\geq m}(M)$ and $\Spinc_{\leq m}(M)$ be the subsets of $\Spinc(M)$ composed of spin$^\c$ structures $\sfr$ with $d(\sfr)=m$, $d(\sfr)\geq m$ and $d(\sfr)\leq m$, respectively.

As shown by Wu in \cite{wu52} and Hirzebruch--Hopf in \cite{hirzebruchhopf58} (before the advent of Seiberg--Witten theory), all spin$^\c$ structures with a zero-dimensional Seiberg--Witten moduli space are induced by almost complex structures, which gives a surjection (see \cite[p. 421]{scorpan05}):
\begin{equation}
\ind\colon
\Jcal(M)
\twoheadrightarrow\Spinc_0(M),
J\mapsto\sfr_J.
\label{eq:ind}
\end{equation}
For $m=0$, Theorem \ref{thrm:main_theorem} combined with this result gives analogous surjections:
\begin{equation*}
p^*\circ\ind\colon
\Jcal\left(M\#(n+3)\oCPtwo\right)
\twoheadrightarrow\Spinc_{\geq 0}(M);
\end{equation*}
\begin{equation*}
q^*\circ\ind\colon
\Jcal\left(M\#(2n+4)\CPtwo\right)\twoheadrightarrow\Spinc_{\leq 2n+4}(M).
\end{equation*}
In the former, increasing the number of blow-ups has no effect on the target, which concerns exactly the spin$^\c$ structures relevant in Seiberg--Witten theory. For $n=0$, Theorem \ref{thrm:main_theorem} further reduces to its strongest special case: 
\begin{crl}
The canonical projection $p\colon M\#3\oCPtwo\twoheadrightarrow M$ induces a surjection:
\begin{equation}
p^*\circ\ind\colon
\Jcal\left(M\#3\oCPtwo\right)
\twoheadrightarrow\Spinc_{\geq 0}(M).
\label{eq:three_blow_ups}
\end{equation}
\label{crl:three_blow_ups}
\end{crl}
First, its implication for the simple type conjecture is:
\begin{thrm}
Every basic class of an oriented closed smooth $4$-manifold with $-b_1+b_2^+$ odd is induced by an almost complex structure on the triple blow-up.
\label{thrm:simple_type}
\end{thrm}
Second, its adaption with the orientation of two blow-ups changed is:
\begin{lem}
The canonical projection $p\colon M\#2\CPtwo\#\oCPtwo\twoheadrightarrow M$ induces surjections:
\begin{equation}
(p^*,d)\colon
\Spinc\left(M\#2\CPtwo\#\oCPtwo\right)
\twoheadrightarrow\Spinc(M)\times 2\Zbb;
\end{equation}
\begin{equation}
p^*\circ\ind\colon
\Jcal\left(M\#2\CPtwo\#\oCPtwo\right)
\twoheadrightarrow\Spinc(M).
\end{equation}
\label{lem:three_blow_ups_twist}
\end{lem}
Combining Lemma \ref{lem:three_blow_ups_twist} with Kotschick's result yields:
\begin{thrm}
For every finitely presented group $G$ there exists a closed almost complex $4$-manifold $M$ with $\pi_1(M)\cong G$ so that $d\colon\Spinc(M)\twoheadrightarrow 2\Zbb$ is surjective.
\label{thrm:combination_with_kotschick}
\end{thrm}

\subsection*{Acknowledgments}

Thanks to Dieter Kotschick for recommending \cite{freeduhlenbeck91} and \cite{donaldsonkronheimer97}.

\section{Setting}
\label{section:setting}

For an oriented closed smooth $4$-manifold $M$, its intersection form canonically pairs the cup product of two second cohomology classes with the fundamental class $[M]\in H_4(M,\Zbb)\cong\Zbb$ (see \parencites[p. 20]{freeduhlenbeck91}[p. 2]{donaldsonkronheimer97}):
\begin{equation*}
Q_M\colon
H^2(M,\Zbb)\times H^2(M,\Zbb)\rightarrow\Zbb,
(x,y)
\mapsto xy[M]
:=\langle x\smile y,[M]\rangle.
\end{equation*}
Let $b_k(M):=\operatorname{rank}H^k(M,\Zbb)\in\Nbb$ be the Betti numbers with Poincaré duality assuring $b_k(M)=b_{4-k}(M)$. Since the intersection form vanishes for all torsion elements, it reduces to a bilinear form on the free part $H^2(M,\Zbb)/\mathrm{tor}\cong\Zbb^{b_2(M)}$. It can therefore be expressed as a symmetric matrix $Q_M\in\Zbb^{b_2(M)\times b_2(M)}$ with $\det(Q_M)=\pm 1$. Let $b_2^\pm(M)$ denote the rank of maximally positive and maximally negative subspaces, respectively. One then has a splitting $b_2(M)=(b_2^++b_2^-)(M)$. Two important topological invariants are the signature and Euler characteristic (see \parencites[p. 22]{freeduhlenbeck91}[p. 3]{donaldsonkronheimer97}), given respectively by:
\begin{equation}
\sigma(M)
:=(b_2^+-b_2^-)(M)
\in\Zbb;
\label{eq:signature}
\end{equation}
\begin{equation}
\chi(M)
:=2-2b_1(M)+b_2(M)
\in\Zbb.
\label{eq:euler_characteristic}
\end{equation}
One has $Q_{\CPtwo}=[+1]$ and $Q_{\oCPtwo}=[-1]$ (see \parencites[p. 21]{freeduhlenbeck91}[p. 3]{donaldsonkronheimer97}). Therefore $\sigma(\CPtwo)=+1$ and $\sigma(\oCPtwo)=-1$. One also has $\chi(\CPtwo)=\chi(\oCPtwo)=3$.

An almost complex structure $J\in\Jcal(M)$ is composed of a complex plane bundle $E\twoheadrightarrow M$ and a real vector bundle isomorphism $E\vert_\mathbb{R}\cong TM$. Its first Chern class $c_1(J):=c_1(E)\in H^2(M,\Zbb)$ fulfills $c_1(J)\operatorname{mod}2=w_2(M)$ and (see \parencites[Thrm. 1.4.15]{gompfstipsicz99}[p. 377]{scorpan05}):
\begin{equation}
c_1^2(J)[M]
=(3\sigma+2\chi)(M).
\label{eq:whh_formula}
\end{equation}
Any such class is induced by an almost complex structure.

A spin$^\c$ structure $\sfr\in\Spinc(M)$ is composed of two complex plane bundles $W^\pm\twoheadrightarrow M$ with $\det(W^+)=\det(W^-)$ and a complex vector bundle isomorphism $TM\otimes\mathbb{C}\cong\underline{\operatorname{Hom}}(W^+,W^-)$ satisfying the Clifford condition. Its first Chern class $c_1(\sfr):=c_1(W^\pm)\in H^2(M,\Zbb)$ fulfills $c_1(\sfr)\operatorname{mod}2=w_2(M)$. Any such class is induced by a spin$^\c$ structure, which is always unique if $H^2(M,\Zbb)$ has no $2$-torsion.

Every almost complex structure $J$ induces a spin$^\c$ structure $\sfr_J$, which establishes the map (\ref{eq:ind}). One has $c_1(J)=c_1(\sfr_J)$, which serves as a unique description of (\ref{eq:ind}) if $H^2(M,\Zbb)$ has no $2$-torsion. However, not every spin$^\c$ structure $\sfr$ is induced by an almost complex structure. Its obstruction comes from (\ref{eq:whh_formula}) and is also reflected in the virtual dimension of its Seiberg--Witten moduli space \cite[Eq. (2.5)]{witten94} (see also \parencites[Thrm. 2.4.24]{gompfstipsicz99}[p. 397]{scorpan05}):
\begin{equation}
d(\sfr)
=\frac{1}{4}\left(c_1^2(\sfr)[M]-(3\sigma+2\chi)(M)\right)
\in 2\Zbb.
\label{eq:sw_moduli_space_dim}
\end{equation}
It is even as van der Blij's lemma (see \cite[p. 170]{scorpan05}) assures $c_1^2(\sfr)[M]-\sigma(M)\in 8\Zbb$ and (\ref{eq:signature}) and (\ref{eq:euler_characteristic}) assure $(\sigma+\chi)(M)=2+2(-b_1+b_2^+)(M)\in 2\Zbb$, respectively.
\begin{rem}
One has (see \parencites[Ex. 1.3.5]{gompfstipsicz99}[p. 118, 120]{scorpan05}):
\begin{equation}
Q_{M\#n\oCPtwo}
\cong Q_M\oplus n[-1];
\label{eq:intersection_blow_up_1}
\end{equation}
\begin{equation}
Q_{M\#n\CPtwo}
\cong Q_M\oplus n[+1];
\label{eq:intersection_blow_up_2}
\end{equation}
\begin{equation}
(3\sigma+2\chi)\left(M\#n\oCPtwo\right)
=(3\sigma+2\chi)(M)-n;
\label{eq:one_blow_up_1}
\end{equation}
\begin{equation}
(3\sigma+2\chi)\left(M\#n\CPtwo\right)
=(3\sigma+2\chi)(M)+5n.
\label{eq:one_blow_up_2}
\end{equation}
\end{rem}
Let $x\in H^2(\CPtwo,\Zbb)\cong H^2(\oCPtwo,\Zbb)\cong\Zbb$ be a generator, meaning $x^2[\CPtwo]=1$ and $x^2[\oCPtwo]=-1$. For $a\in\Zbb$, let $\osfr_a\in\Spinc(\oCPtwo)$ be the spin$^\c$ structure with $c_1(\osfr_a)=(2a+1)x$ and let $\sfr_a\in\Spinc(\CPtwo)$ be the spin$^\c$ structure with $c_1(\sfr_a)=(2a+1)x$. One has:
\begin{equation}
c_1^2(\osfr_a)\left[\oCPtwo\right]
=-(2a+1)^2
=-4a(a+1)-1;
\label{eq:intersection_form_cp2}
\end{equation}
\begin{equation}
c_1^2(\sfr_a)[\CPtwo]
=(2a+1)^2
=4a(a+1)+1.
\label{eq:intersection_form_cp2_bar}
\end{equation}
Let $\osfr_{-1,0}\in\{\osfr_{-1},\osfr_0\}$, $\sfr_{-1,0}\in\{\sfr_{-1},\sfr_0\}$ and $\sfr_{-2,1}\in\{\sfr_{-2},\sfr_1\}$, then one has $c_1^2(\osfr_{-1,0})[\oCPtwo]=-1$, $c_1^2(\sfr_{-1,0})[\CPtwo]=1$ and $c_1^2(\sfr_{-2,1})[\CPtwo]=9$, respectively.
\begin{lem}
For an almost complex $4$-manifold $M$:
\begin{itemize}
\item[a)] $M\#\oCPtwo$ is almost complex.
\item[b)] $M\#2\CPtwo$ is almost complex.
\end{itemize}
\label{lem:almost_complex}
\end{lem}
\begin{proof}
Let $J\in\Jcal(M)$. Using (\ref{eq:whh_formula}), (\ref{eq:intersection_blow_up_1}), (\ref{eq:one_blow_up_1}) and (\ref{eq:intersection_form_cp2}):
\begin{equation}
(c_1(J)\#\osfr_{-1,0})^2\left[M\#\oCPtwo\right]
=(3\sigma+2\chi)(M\#\oCPtwo).
\end{equation}
Analogous from (\ref{eq:whh_formula}), (\ref{eq:intersection_blow_up_2}), (\ref{eq:one_blow_up_2}) and (\ref{eq:intersection_form_cp2_bar}):
\begin{equation}
(c_1(J)\#c_1(\sfr_{-1,0})\# c_1(\sfr_{-2,1}))^2[M\#2\CPtwo]
=(3\sigma+2\chi)(M\#2\CPtwo).
\end{equation}
Hence $c_1(J)\#c_1(\osfr_{-1,0})$ and $c_1(J)\# c_1(\sfr_{-1,0})\# c_1(\sfr_{-2,1})$ fulfill (\ref{eq:whh_formula}). According to (\ref{eq:ind}), both come from an almost complex structure.
\end{proof}

\section{Proofs}
\label{section:proofs}

Let $M$ be an oriented smooth $4$-manifold and $\sfr\in\Spinc(M)$. Using (\ref{eq:sw_moduli_space_dim}), (\ref{eq:intersection_blow_up_1}), (\ref{eq:one_blow_up_1}) and (\ref{eq:intersection_form_cp2}) (see \cite[p. 459]{nicolaescu00}):
\begin{equation}
d(\sfr\#\osfr_a)
=d(\sfr)
-a(a+1).
\label{eq:sw_moduli_space_dim_blow_up_1}
\end{equation}
Analogous from (\ref{eq:sw_moduli_space_dim}), (\ref{eq:intersection_blow_up_2}), (\ref{eq:one_blow_up_2}) and (\ref{eq:intersection_form_cp2_bar}):
\begin{equation}
d(\sfr\#\sfr_a)
=d(\sfr)
+a(a+1)-1.
\label{eq:sw_moduli_space_dim_blow_up_2}
\end{equation}
Relying on the virtual dimension of all Seiberg--Witten moduli spaces being even, the core of the proofs is the combination of (\ref{eq:sw_moduli_space_dim_blow_up_1}) and (\ref{eq:sw_moduli_space_dim_blow_up_2}) with the following result from number theory:
\begin{thrm}
[Gauss' Eureka theorem]
For every $n\in\Nbb$, there exist $a,b,c\in\Nbb$ with:
\begin{equation}
2n
=a(a+1)
+b(b+1)
+c(c+1).
\label{eq:gauss_eureka_thrm}
\end{equation}
\label{thrm:gauss_eureka_thrm}
\end{thrm}
\begin{proof}
[Proof of Theorem \ref{thrm:main_theorem}]
Let $\sfr\in\Spinc(M)$ and $m\in\Zbb$ fulfill $2m\leq d(\sfr)$. According to Theorem \ref{thrm:gauss_eureka_thrm}, there exist $a,b,c\in\Nbb$ such that:
\begin{equation}
d(\sfr)
-2m
=a(a+1)
+b(b+1)
+c(c+1).
\label{eq:sw_moduli_space_eureka_1}
\end{equation}
For $\widehat\sfr=\sfr\#\osfr_a\#\osfr_b\#\osfr_c\in\Spinc(M\#3\oCPtwo)$, applying (\ref{eq:sw_moduli_space_dim_blow_up_1}) to each of the three added spin$^\c$ structures makes it a preimage:
\begin{equation*}
(p^*,d)(\widehat\sfr)
=(\sfr,2m).
\end{equation*}
Hence the map (\ref{eq:main_theorem_1}) is surjective. It is also well-defined and furthermore both for $n\geq 0$ since (\ref{eq:sw_moduli_space_dim_blow_up_1}) with $a(a+1)\geq 1$ and the minimum attained for $a=-1,0$ yields:
\begin{equation*}
d(\widehat\sfr\#n\osfr_{-1,0})
=d(\widehat\sfr)
\leq d(\sfr).
\end{equation*}
Let $m\in\Zbb$ fulfill $d(\sfr)\leq 2m+4$. According to Theorem \ref{thrm:gauss_eureka_thrm}, there exist $a,b,c\in\Nbb$ such that:
\begin{equation}
2m+4-d(\sfr)
=a(a+1)+b(b+1)+c(c+1).
\label{eq:sw_moduli_space_eureka_2}
\end{equation}
For $\widehat\sfr=\sfr\#\sfr_a\#\sfr_b\#\sfr_c\#\sfr_{-1,0}\in\Spinc(M\# 4\CPtwo)$, applying (\ref{eq:sw_moduli_space_dim_blow_up_2}) to each of the four added spin$^\c$ structures makes it a preimage:
\begin{equation*}
(q^*,d)(\widehat\sfr)
=(\sfr,2m).
\end{equation*}
Hence the map (\ref{eq:main_theorem_2}) is surjective. It is also well-defined and furthermore both for $n\geq 0$ since (\ref{eq:sw_moduli_space_dim_blow_up_2}) with $a(a+1)\geq 1$ and the minimum attained for $a=-1,0$ yield:
\begin{equation*}
d(\widehat\sfr\#2n\sfr_{-1,0})
+2n
\geq d(\widehat\sfr)
\geq d(\sfr)-4.
\end{equation*}
\end{proof}
\begin{proof}
[Proof of Lemma \ref{lem:three_blow_ups_twist}]
If the sign of the entire last summand in (\ref{eq:gauss_eureka_thrm}) is changed, then it directly holds for every $n\in\Zbb$ using $a=-1,0$, $b=n$ and $c=n-1$. Let $\sfr\in\Spinc(M)$, $m\in\Zbb$ and:
\begin{equation*}
n=-\frac{d(\sfr)}{2}+m+1.
\end{equation*}
For $\widehat\sfr=\sfr\#\sfr_{-1,0}\#\sfr_n\#\osfr_{n-1}\in\Spinc\left(M\#2\CPtwo\#\oCPtwo\right)$, applying (\ref{eq:sw_moduli_space_dim_blow_up_1}) once and (\ref{eq:sw_moduli_space_dim_blow_up_2}) twice makes it a preimage:
\begin{equation*}
(p^*,d)(\widehat\sfr)
=(\sfr,2m).
\end{equation*}
\end{proof}
\begin{proof}
[Proof of Theorem \ref{thrm:combination_with_kotschick}]
Let $G$ be a finitely presented group. According to \cite{kotschick92}, there exists an almost complex closed $4$-manifold $M$ with $\pi_1(M)\cong G$. According to Lemma \ref{lem:almost_complex} (or implied by Lemma \ref{lem:three_blow_ups_twist}), $M\#2\CPtwo\#\oCPtwo$ is also an almost complex closed $4$-manifold. Furthermore $\pi_1(M\#2\CPtwo\#\oCPtwo)\cong G$ since $\CPtwo$ and $\oCPtwo$ are simply connected. Lemma \ref{lem:three_blow_ups_twist} then yields the claim.
\end{proof}

\section{Discussion}
\label{section:discussion}

John A. Ewell provided two results related to Theorem \ref{thrm:gauss_eureka_thrm}: In \cite{ewell92rep}, he gives a criterium for special cases reducing it to two summands, which transfers directly to reducing Corollary \ref{crl:three_blow_ups} to two blow-ups. However, this requires considering prime factors, which for arbitrary differences like in (\ref{eq:sw_moduli_space_eureka_1}) and (\ref{eq:sw_moduli_space_eureka_2}) are completely unknown. In \cite{ewell92sum}, he gives an explicit formula for how many possible decompositions there are, which transfers directly to how much Corollary \ref{crl:three_blow_ups} deviates from injectivity.
\begin{rem}
Three blow-ups in Corollary \ref{crl:three_blow_ups} is optimal in general. For instance, for $n=5$, one has $a=1$, $b=1$ and $c=2$ in (\ref{eq:gauss_eureka_thrm}), making it the smallest example of
Theorem \ref{thrm:gauss_eureka_thrm} for which no $a,b,c\in\Nbb$ can vanish. For $M=\CPtwo$, using (\ref{eq:sw_moduli_space_dim}) and (\ref{eq:intersection_form_cp2_bar}), one has $d(\sfr_{-4,3})=10$ (see also \cite[p. 168]{nicolaescu00}). From (\ref{eq:sw_moduli_space_dim}), (\ref{eq:intersection_blow_up_2}) and (\ref{eq:one_blow_up_1}) one has:
\begin{align*}
d(\sfr_{-4,3}\#\sfr_a\#\sfr_b)
&=\frac{1}{4}\left(
c_1^2(\sfr_{-4,3}\#\sfr_a\#\sfr_b)\left[\CPtwo\#2\oCPtwo\right]
-(3\sigma+2\chi)\left(\CPtwo\#2\oCPtwo\right)\right) \\
&=10-a(a+1)-b(b+1)
\neq 0.
\end{align*}
Hence no $a,b\in\Nbb$ exist so that $\sfr_{-4,3}\#\sfr_a\#\sfr_b\in\Spinc(\CPtwo\#2\oCPtwo)$ fulfills (\ref{eq:whh_formula}) and is therefore induced by an almost complex structure according to (\ref{eq:ind}).
\end{rem}
It is also worth noting that both results have been used in other fields using the same concepts as this paper: For instance, they were used to study circle actions on compact almost complex manifolds by Godinho--Pelayo--Sabatini in \cite{godinhopelayosabatini14}, almost complex closed manifolds with the same rational cohomology as a connected sum of $\Hbb P^2$ and $\overline{\Hbb P^2}$ by Milivojević in \cite{mil22}, generalizations of spin$^\c$ structures by Albanese--Milivojević in \cite{albmil22}, and string$^\c$ structures by Duan--Han--Huang in \cite{duanhanhuang20}, among others.

Two further remarks are important concerning Theorem \ref{thrm:simple_type}. First, the simple type conjecture requires $b_2^+\geq 2$, which cannot be omitted. Otherwise, $\CPtwo$ provides a counterexample since it fulfills all other conditions, but has $b_2^+(\CPtwo)=1$ and basic classes $\sfr$ with $d(\sfr)>0$ (see  \cite[p. 168]{nicolaescu00}). However, Theorem \ref{thrm:simple_type} does omit $b_2^+\geq 2$. Second, blow-ups with the spin$^\c$ structures $\sfr_{-1,0}$ preserve the Seiberg--Witten invariant while all others cause it to vanish (see \cite[Thrm. 2.4.9]{gompfstipsicz99}). However, in general these are considered in the proofs. Both of these aspects are based on the fact that Seiberg--Witten invariants, which require $b_2^+\geq 2$ for metric independence and which are at the core of the simple type conjecture, were not used to derive Theorem \ref{thrm:simple_type}. An immediate question is then whether Theorem \ref{thrm:simple_type} can be improved by using the Seiberg--Witten invariants and Ewell's results to obtain an even closer approximation to the simple type conjecture:
\begin{con}
Every basic class of an oriented closed smooth $4$-manifold with $-b_1+b_2^+$ odd and $b_2^+\geq 2$ is induced by an almost complex structure on the double blow-up.
\end{con}

\printbibheading
\printbibliography[heading=none]

\end{document}